\documentclass[preprint,12pt]{elsarticle}

\usepackage{amssymb}
\usepackage{amsmath,amssymb,amsfonts}
\usepackage{tikz}
\usepackage{tabularx}

\newtheorem{theorem}{Theorem}
\newtheorem{corollary}{Corollary}
\newtheorem{lemma}{Lemma}
\newdefinition{rmk}{Definition}
\newproof{proof}{Proof}
\newtheorem{remark}{Remark}

\journal{arXiv}

\begin{document}

\begin{frontmatter}


\title{Exponential stability of Euler-Bernoulli beam under boundary controls in rotation and angular velocity}

\author[1]{Alemdar Hasanov\fnref{fn1}}
\ead{alemdar.hasanoglu@gmail.com}

\cortext[cor1]{Corresponding author}
\fntext[fn1]{Department of Mathematics, Kocaeli University, Turkey}

\address[1]{Department of Mathematics, Kocaeli University, Turkey}

\begin{abstract}
This paper addresses the analysis of a boundary feedback system involving a non-homogeneous Euler-Bernoulli beam governed by the equation $m(x)u_{tt}+\mu(x)u_{t}$$+\left(r(x)u_{xx}\right)_{xx}=0$, subject to the
initial $u(x,0)=u_0(x)$, $u_t(x,0)=v_0(x)$ and boundary conditions $u(0,t)=0$, $\left (-r(x)u_{xx}(x,t)\right )_{x=0}=-k^{-}_r u_{x}(0,t)-k^{-}_a u_{xt}(0,t)$,
$u(\ell,t)=0$, $\left (-r(x)u_{xx}(x,t)\right )_{x=\ell}=-k^{+}_r u_{x}(\ell,t)-k^{+}_a u_{xt}(\ell,t)$, with boundary control at both ends resulting from the rotation and angular velocity. The approach proposed in this study relies on the utilization of regular weak solutions, energy identity, and a physically motivated Lyapunov function. By imposing natural assumptions concerning physical parameters and other inputs, which ensure the existence of a regular weak solution, we successfully derive a uniform exponential decay estimate for the system's energy. The decay rate constant featured in this estimate is solely dependent on the physical and geometric properties of the beam. These properties encompass crucial parameters such as the viscous external damping coefficient $\mu(x)$, as well as the boundary springs $k^{-}_r,k^+_r $ and dampers $k^{-}_a,k^+_a$. To illustrate the practical effectiveness of our theoretical findings, numerical examples are provided. These examples serve to demonstrate the applicability and relevance of our derived results in real-world scenarios.
\end{abstract}


\end{frontmatter}


\section{Introduction}
\label{intro}

Submarine pipelines and long bridges can be considered as an elastic beam with both ends controlled by the boundary rotation and angular velocity \cite{Cai:2022, Liu:2018}. In many studies related to pipeline modeling, the pipes are defined as beams resting on a rigid seabed without any penetration (see \cite{Hong:2015} and references therein). However, such hypotheses are not always satisfied in practice. An analysis of the torsional effects on pipe lateral buckling was given in \cite{Grognec:2020}, where essential influence of torsion under some specific boundary conditions was demonstrated analytically. Similar situation arise in bridge models governed by the Euler-Bernoulli beam. Namely, it is very important for the sensitivity analysis of bridges to obtain a relationship between the rotation spring constant and the bridge responses (deflections/slopes). This relationship can then be used for evaluating the support condition of bridges \cite{Park:2019}. Furthermore, in modeling of long flexible structures through the Euler-Bernoulli equation, the bending moment at the end of the beam is controlled by the linear feedback of rotation angle and angular velocity, and the shear force at the same end is controlled by the linear feedback of displacement
and velocity. We refer \cite{F-Guo:2004} and references therein, for the detailed description of such models.

Considering the effect of the above factor on both models, there is a need for a realistic model that will take into account the effects of both the rotation spring and the angular velocity damper at both ends of the beam, within the framework of the Euler-Bernoulli beam equation. In the most natural way, this can be taken into account by the corresponding boundary conditions at both ends of the beam, including a linear combinations of the rotation spring and the angular velocity damper. This leads to the following mathematical  model:
\begin{eqnarray}\label{1}
\left\{ \begin{array}{ll}
m(x)u_{tt}+\mu(x)u_{t}+\left(r(x)u_{xx}\right)_{xx}=0,\, (x,t) \in \Omega_T, \\ [4pt]
    u(x,0)=u_0(x), \, u_t(x,0)=u_1(x),\, x\in (0,\ell),\\ [4pt]
	u(0,t)=0,\, \left (-r(x)u_{xx}(x,t)\right )_{x=0}=-k^{-}_r u_{x}(0,t)-k^{-}_a u_{xt}(0,t), \\ [4pt]
\quad u(\ell,t)=0,\, \left (-r(x)u_{xx}(x,t)\right )_{x=\ell}=k^{+}_r u_{x}(\ell,t)+k^{+}_a u_{xt}(\ell,t),\\ [4pt]
\qquad \qquad  \qquad \qquad  \qquad \qquad \qquad \qquad \qquad  \qquad \qquad ~t\in [0,T],
\end{array} \right.
\end{eqnarray}
where $\Omega_T=(0,\ell)\times(0,T)$, $\ell>0$ is the length of the beam and $T>0$ is the final time.

Here and below, $u(x,t)$ is the deflection, $u_t(x,t)$, $u_x(x,t)$, $u_{xt}(x,t)$, $u_{xx}(x,t)$, $-\left(r(x)u_{xx}\right)$ and $-\left(r(x)u_{xx}\right)_{x}$ are the velocity, rotation (or slope), angular velocity, curvature, moment and shear force, respectively \cite{Clough-Penzien:1975, Inman:2015}. Further, $m(x)=\rho(x)S(x)>0$, while $\rho(x)$ is the mass density and $S(x)$ is the cross section area of the beam, and $r(x):=E(x)I(x)>0$ represent the flexural rigidity (or bending stiffness) of the beam, respectively, while $E(x)>0$ is the elasticity modulus and $I(x)>0$ is the moment of inertia. The non-negative coefficient $\mu(x):=\gamma\, m(x)$ of viscous resistance to transverse motion of the beam represents the viscous external damping, while $\gamma \ge 0$ is the damping constant of proportionality \cite{Banks:Inman:1991}. Furthermore, nonnegative constants $k^{-}_r, k^{-}_a \ge 0$ and $k^{+}_r, k^{+}_a \ge 0$ are the stiffness of the torsional springs and dampers on the left and right ends of the beam, respectively.

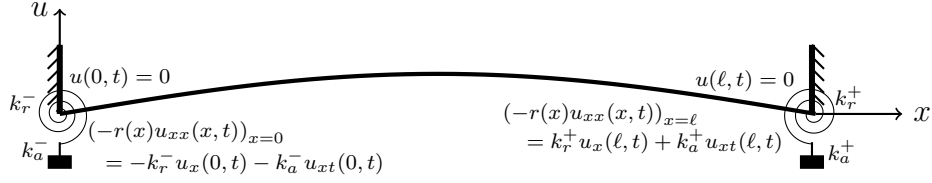
\begin{figure}
  \centering
 \hspace*{0.6cm} \begin{tikzpicture}[scale=.8]
      \draw[line width=0.8mm, black] (0.0,0.0) -- (0.0,1.15);
      \draw[color=black,thick] (0.0,1.15) -- (-0.20,0.95);
      \draw[color=black,thick] (0.0,0.95) -- (-0.20,0.75);
      \draw[color=black,thick] (0.0,0.75) -- (-0.20,0.55);
      \draw[color=black,thick] (0.0,0.55) -- (-0.20,0.35);
      \draw[color=black,thick] (0.0,0.35) -- (-0.20,0.15);
      \draw[line width=0.8mm, black] (12.5,0.0) -- (12.5,1.15);
          \draw[color=black,thick] (12.5,1.15) -- (12.7,0.95);
          \draw[color=black,thick] (12.5,0.95) -- (12.7,0.75);
          \draw[color=black,thick] (12.5,0.75) -- (12.7,0.55);
          \draw[color=black,thick] (12.5,0.55) -- (12.7,0.35);
          \draw[color=black,thick] (12.5,0.35) -- (12.7,0.15);
      \draw[thick,->] (4*3.14,0) -- (14,0) node [right]{$x$};
      \draw[thick,->] (0,0) -- (0,1.75) node [left]{$u$};
      \draw[color=black,ultra thick,smooth,domain=0:{4*3.14}] plot (\x,{sin(0.25*deg(\x))/1.5});
      \node[label=right:{\scriptsize$ u(0,t)=0$}] at (-0.2,0.6) {};
      \node[label=right:{\scriptsize$\left (-r(x)u_{xx}(x,t)\right )_{x=0}$}] at (0.1,-0.3) {};
      \node[label=right:{\scriptsize$=-k^{-}_r u_{x}(0,t)-k^{-}_a u_{xt}(0,t)$}] at (0.4,-0.8) {};
      \node[label=right:{\scriptsize$ u(\ell,t)=0$}] at (10.2,0.55) {};
       \node[label=right:{\scriptsize$\left (-r(x)u_{xx}(x,t)\right )_{x=\ell}$}] at (7.0,0.0) {};
      \node[label=right:{\scriptsize$=k^{+}_r u_{x}(\ell,t)+k^{+}_a u_{xt}(\ell,t)$}] at (7.4,-0.45) {};
       \node[label=right:{\scriptsize$k^{-}_r$}] at (-1.2,0.2) {};
       \node[label=right:{\scriptsize$k^{-}_a$}] at (-1.0,-0.58) {};
        \node[label=right:{\scriptsize$k^{+}_r$}] at (12.52,0.28) {};
       \node[label=right:{\scriptsize$k^{+}_a$}] at (12.4,-0.63) {};
      \def\spiral[#1](#2)(#3:#4:#5){
\pgfmathsetmacro{\domain}{pi*#3/180+#4*2*pi}
\draw [#1,
       shift={(#2)},
       domain=0:\domain,
       variable=\t,
       smooth,
       samples=int(\domain/0.08)] plot ({\t r}: {#5*\t/\domain})}
       \draw [samples=100,smooth,domain=10:26.7] 
       plot ({-\x r}:{0.0007*\x*\x});      
  \def\x{27.5}
       \draw[thick] (0.0,-0.47)-- (0.0,-0.75);
       \draw[line width=1.8mm, black] (-0.2,-0.8)-- (0.2,-0.8);
       \spiral[black](12.5,0)(-1.5:2.76:0.5);
       \draw[thick] (12.5,-0.47)-- (12.5,-0.75);
       \draw[line width=1.8mm, black] (12.3,-0.8)-- (12.7,-0.8);
  \end{tikzpicture}
  \caption{Beam connected to torsional springs and dampers at both ends}
\label{Fig-1}
  \end{figure}

The boundary conditions $\left (-r(x)u_{xx}(x,t)\right )_{x=0}=-k^{-}_r u_{x}(0,t)-k^{-}_a u_{xt}(0,t)$ and $\left (-r(x)u_{xx}(x,t)\right )_{x=\ell}=k^{+}_r u_{x}(\ell,t)+k^{+}_a u_{xt}(\ell,t)$ at the left and right ends of the beam, respectively, mean the controls resulting from the linear combination of rotation and angular velocity. In this context, the above parameters $k^{-}_r,\, k^{-}_a,\,k^{+}_r,\,k^{+}_a$ are defined also as the boundary controls.

Geometry of the problem (\ref{Fig-1}) is given in Fig. \ref{Fig-1}.

This work is devoted to the systematic study of the following issues. \emph{Under what minimum conditions imposed on the input data is the energy of the system governed by (\ref{1}) exponentially stable?} \emph{If the system governed by (\ref{1}) is stable, how much does each damping parameter $\gamma$,  $k^{-}_a$ and $k^{+}_a$ contribute to this stability?} It should be especially noted that the nature of both the external and the boundary damping mechanisms greatly changes the nature of the vibration, and hence controls the response of the beam, as the experimental and theoretical results discussed in  \cite{Banks:Inman:1991, Crandall:1970} show.

Modeling of large flexible structures through a class of Euler-Bernoulli beams with structural damping, has begun to be developed, starting with studies \cite{Chen-Krantz:1988, Chen-Russell:1982, Russell:1978}. The exponential stability of distributed systems governed by Euler-Bernoulli beam equation under classical boundary conditions has been discussed starting from the work  \cite{Chen-Russell:1982}, and then more general results are obtained in \cite{Chen-Krantz:1988, Huang:1986, Huang:1988}. Various methods have been developed in the literature for initial boundary value problems for Euler-Bernoulli equations with a boundary feedback systems. Among these methods, the spectral method turned out to be efficient and useful since it allows to establish the Riesz basis property, which is the most fundamental property of a linear vibrating system \cite{Chen-Xu:2014, Guo:2001, Guo:2002, F-Guo:2004}. In turn, this property means that the generalized eigenvectors of the system form an unconditional basis of the (state) Hilbert space. With semigroup approach, this allows to derive the spectrum determined growth condition and the exponential stability for a system.

In the exponential stability estimate $\mathcal{E}(t) \le M e^{-\omega t} \mathcal{E}(0)$ obtained in the studies listed above, the relationship of the decay rate parameter $\omega >0$ with the physical and geometric parameters of the beam, including the damping coefficient $\mu(x) \ge 0$ and the stiffness $k^{-}_a, k^{+}_a \ge 0$ of the torsional dampers, has not been determined. Since the relationship of this decay rate parameter with the damping parameters is not known, in concrete applications, such an evaluation does not give a qualified result.
In this paper, we develop the approach based on the weak solution theory for the initial boundary value problem (\ref{1}), energy estimates and the Lyapunov method to establish an exponential stability estimate for system (\ref{1}) under minimum conditions imposed on the input data. Furthermore, this approach allows us to derive the role of both types of parameters in the exponential decay of the solution. To our knowledge, this model, defined by the initial boundary value problem (\ref{1}), in which the viscous external and boundary (torsional) damping factors are considered together and in the presence of torsional springs, is discussed for the first time in the literature.

The rest of the paper is structured as follows. Energy identity and dissipativity of system (\ref{1}) are derived in Section 2. In Section 3, the Lyapunov function is introduced and then energy decay estimate for system (\ref{1}) is derived. Numerical examples are presented in Section 4. Some concluding remarks are given in the final Section 5.

\section{Necessary estimates for the weak solution of problem (\ref{1})}

We assume that the inputs in (\ref{1}) satisfy the following basic conditions:
\begin{eqnarray} \label{2}
\left \{ \begin{array}{ll}
\rho_S, \mu,r \in L^\infty(0,\ell),\\ [3pt]
0<m_0\leq m(x)\leq m_1,~0\leq \mu_0\leq \mu(x)\leq \mu_1,\\ [3pt]
0<r_0\leq r(x)\leq r_1,\,  x\in (0,\ell),\\ [3pt]
u_0\in  H^2(0,\ell),~u_1\in  L^2(0,\ell),\\ [3pt]
k^{-}_r, k^{-}_a,k^{+}_r, k^{+}_a \ge 0,\\ [3pt]
\gamma+k^{-}_r+ k^{-}_a+k^{+}_r+k^{+}_a >0.
\end{array} \right.
\end{eqnarray}

For the case when all the parameters $k^{-}_r, k^{-}_a,k^{+}_r, k^{+}_a$ are equal to zero, under conditions (\ref{2}), the existence of the weak solution $u\in L^2(0,T; \mathcal{V}^2(0,\ell))$, with $u_t\in L^2(0,T;L^2(0,\ell))$ and $u_{tt}\in L^2(0,T;H^{-2}(0,\ell))$ of the initial boundary value problem (\ref{1}) was proved in   \cite{Hasanov-Romanov:2021}. Here and below,
\begin{eqnarray*}
\mathcal{V}^2(0,\ell):=\{v\in H^2(0,\ell):\, v(0)=v(\ell)=0,\},
\end{eqnarray*}
and $H^2(0,\ell)$ is the Sobolev space \cite{Evans:2002}. For system (\ref{1}), with $k^{-}_r, k^{-}_a,k^{+}_r, k^{+}_a>0$, the existence of the weak solution $u\in L^2(0,T; \mathcal{V}^2(0,\ell))$ can be proved in the similar way. In this section we derive necessary energy identities and estimates for the weak solution of problem (\ref{1}).

\medskip
\begin{theorem}\label{Theorem-1}
Assume that the inputs in (\ref{1}) satisfy the basic conditions (\ref{2}).
Then the following energy identity holds:
\begin{eqnarray}\label{3}
\mathcal{E}(t) + \int_0^t \int_0^\ell \mu(x) u_{\tau}^2 (x,\tau) dx d \tau \qquad \qquad \qquad    \qquad \qquad \qquad \qquad \nonumber \\ [1pt]
\qquad \qquad  =\mathcal{E}(0)-k^{-}_a \int_0^t u_{x\tau}^2(0,\tau) d \tau-k^{+}_a \int_0^t u_{x\tau}^2(\ell,\tau) d \tau,\,t\in[0,T],
\end{eqnarray}
where
\begin{eqnarray}\label{4}
\mathcal{E}(t)=\frac{1}{2} \int_0^\ell \left [ m(x) u^2_{t}(x,t)+r(x) u^2_{xx}(x,t)\right ] dx \qquad \qquad \qquad \nonumber \\ [1pt]
+\frac{1}{2}\,k^{-}_r u_{x}^2(0,t) +\frac{1}{2}\,k^{+}_r\, u_{x}^2 (\ell,t),~t\in[0,T],
\end{eqnarray}
is the total energy of system (\ref{1}) and
\begin{eqnarray}\label{5}
\mathcal{E}(0)=\frac{1}{2} \int_0^\ell \left [ m(x)\left ( u_{1}(x)\right)^2 +
r(x) \left ( u''_{0}(x)\right)^2\right ] dx \qquad \qquad \qquad \nonumber \\ [1pt]
\qquad +\frac{1}{2}\, k^{-}_r \left ( u'_{0}(0)\right)^2+\frac{1}{2}\, k^{+}_r \left ( u'_{0}(\ell)\right)^2
\end{eqnarray}
is the initial value of the total energy.
\end{theorem}
{\bf Proof.} Multiply both sides of equation (\ref{1}) by $u_t(x,t)$, integrate it over $\Omega_t:=(0,\ell)\times (0,t)$, employ the identity
\begin{eqnarray}\label{6}
\int_0^t\int_0^\ell (r(x)u_{xx})_{xx} u_{\tau} dx d\tau=  \int_0^t\int_0^\ell [(r(x)u_{xx})_x u_{\tau}-r(x)u_{xx} u_{x\tau}]_x dx d\tau  \nonumber \\ [1pt]
 + \,\frac{1}{2}\int_0^t\int_0^\ell \left (r(x)u_{xx}^2\right )_{\tau} dx d\tau, \quad
\end{eqnarray}
$t \in (0,T]$. Then we obtain the following integral identity:
\begin{eqnarray*}
\frac{1}{2} \int_0^t\int_0^\ell \left (\rho_S(x) u_{\tau}^2\right )_{\tau}dx\,d\tau +\frac{1}{2} \int_0^t \int_0^\ell \left (r(x)u_{xx}^2\right )_{\tau}dx\,d\tau \qquad \qquad \qquad  \\ [1pt]
\qquad +\int_0^t \left ((r(x)u_{xx})_x u_{\tau}-r(x)u_{xx} u_{x\tau} \right)_{x=0}^{x=\ell} d\tau +\int_0^t \int_0^\ell \mu(x) u_{\tau}^2 dx d \tau=0,
\end{eqnarray*}
for all $t \in (0,T]$. Using here the initial and boundary conditions (\ref{1}), we get:
 \begin{eqnarray*}
\frac{1}{2} \int_0^\ell \left [m(x) u^2_{t}+ r(x) u_{xx}\right ]dx
+\frac{1}{2}\,k^{-}_r u_{x}^2(0,t) +\frac{1}{2}\,k^{+}_r\, u_{x}^2 (\ell,t)\qquad \qquad \qquad  \\ [1pt]
+\int_0^t \int_0^\ell \mu(x) u_{\tau}^2 dx d \tau \\ [1pt]
= \frac{1}{2} \int_0^\ell \left [m(x)\left ( u_{1}(x)\right)^2 +
r(x) \left ( u''_{0}(x)\right)^2\right ] dx  +\frac{1}{2}\, k^{-}_r \left ( u'_{0}(0)\right)^2+\frac{1}{2}\, k^{+}_r \left ( u'_{0}(\ell)\right)^2 \\ [1pt]
-k^{-}_a \int_0^t u_{x\tau}^2(0,\tau) d \tau-k^{+}_a \int_0^t u_{x\tau}^2(\ell,\tau) d \tau,\,t\in[0,T],
\end{eqnarray*}
for all $t \in (0,T]$. This leads to (\ref{3}) with (\ref{4}) and (\ref{5}).\hfill$\Box$

\medskip
\begin{remark}\label{Remark-1}
The integral identity (\ref{3}), with (\ref{4}) and (\ref{5}), clearly shows that the increase in the stiffness of the torsional springs $k^{-}_r$ and $k^{+}_r$ leads to an increase in the total energy $\mathcal{E}(t)$. Conversely, the increase in the stiffness of the torsional dampers  $k^{-}_a$ and $k^{+}_a$ leads to a decrease in the total energy.
\end{remark}

\medskip
\begin{remark}\label{Remark-2}
The sum
\begin{eqnarray*}
\frac{1}{2}\,k^{-}_r u_{x}^2(0,t) +\frac{1}{2}\,k^{+}_r\, u_{x}^2 (\ell,t),~t\in[0,T]
\end{eqnarray*}
in (\ref{4}) represents the energy of the rigid motion of the elastic system (\ref{4}), generated by the  spring constants $k^{-}_r,k^{+}_r\ge0$.
\end{remark}

\medskip
\begin{lemma}\label{Lemma-1}
Assume that the basic conditions (\ref{2}) hold. Then for the decay rate of the total energy the following integral formula is valid:
\begin{eqnarray}\label{7}
\frac{d \mathcal{E}(t)}{dt} =-\int_0^\ell \mu(x) u^2_{t}dx- k^{-}_a u_{xt}^2 (0,t) - k^{+}_a u_{xt}^2(\ell,t),\, t\in (0,T).
\end{eqnarray}
\end{lemma}
{\bf Proof.} From formula (\ref{4}) for the total energy we deduce that
\begin{eqnarray*}
\frac{d \mathcal{E}(t)}{dt}= \int_0^\ell \left [ m(x) u_{t}u_{tt}+r(x) u_{xx}u_{xxt}\right ] dx \qquad \qquad \qquad \qquad \nonumber\\ [1pt]
\qquad \qquad  \qquad +k^{-}_r u_{x}(0,t)u_{xt}(0,t) +
k^{+}_r u_{x}(\ell,t)u_{xt}(\ell,t),~t\in[0,T].
\end{eqnarray*}
Using here the identities
\begin{eqnarray*}
\int_0^\ell m(x)u_{t} u_{tt} dx =-\int_0^\ell \mu (x) u_{t}^2 dx
-\int_0^\ell  \left (r(x) u_{xx}\right )_{xx} u_{t} dx,  \qquad \qquad\\ [2pt]
\int_0^\ell  \left (r(x) u_{xx}\right )_{xx} u_{t} dx=
\int_0^\ell  r(x) u_{xx}u_{xxt}dx +k^{-}_r u_{x}(0,t)u_{xt}(0,t)\qquad \\ [1pt]
+k^{-}_a u^2_{xt}(0,t)+k^{+}_r u_{x}(\ell,t)u_{xt}(\ell,t)
+k^{+}_a u^2_{xt}(\ell,t),~t\in[0,T],
\end{eqnarray*}
we arrive at the required result (\ref{7}).
 \hfill$\Box$

\medskip
\begin{corollary}\label{Corollary-1}
Integrating (\ref{7}) over $(0,t)$ we arrive at the energy identity introduced in (\ref{3}), that is
\begin{eqnarray}\label{8}
\mathcal{E}(t) =\mathcal{E}(0)-\int_0^t\int_0^\ell \mu(x) u^2_{\tau}(x,\tau)dx d\tau \qquad \qquad \qquad \qquad \nonumber\\ [1pt]
- \int_0^t \left [k^{-}_a u_{x \tau}^2 (0,\tau) + k^{+}_a u_{x \tau}^2(\ell,t) \right ]d\tau,~t\in [0,T].
\end{eqnarray}
In particular,
\begin{eqnarray*}
\mathcal{E}(t) \le \mathcal{E}(0),~t\in[0,T],
\end{eqnarray*}
that is, the energy of the system (\ref{1}) is dissipating.
\end{corollary}

\section{Lyapunov function and exponential stability estimate}

Introduce the auxiliary function:
\begin{eqnarray}\label{9}
\mathcal{J}(t)= \int_0^\ell m(x) u\,u_{t}dx+\frac{1}{2} \int_0^\ell \mu(x) u^2dx+ \frac{1}{2}\, k^{-}_a  u_x^2 (0,t) +\frac{1}{2}\, k^{+}_a  u_x^2 (\ell,t),
\end{eqnarray}
$t\in[0,T]$, that includes all the damping parameters.

\medskip
\begin{lemma}\label{Lemma-2}
Assume that the basic conditions (\ref{2}) are satisfied. Then between the auxiliary function $\mathcal{J}(t)$ and the energy function $\mathcal{E}(t)$, the following relationship
holds:
\begin{eqnarray}\label{10}
\frac{d \mathcal{J}(t)}{dt}= 2 \int_0^\ell m(x) u_t^2dx -2\mathcal{E}(t),~t\in[0,T].
\end{eqnarray}
\end{lemma}
{\bf Proof.} Taking the derivative of the function $\mathcal{J}(t)$ with respect to the time variable and using then the equation (\ref{1}) we find:
\begin{eqnarray*}
\frac{d \mathcal{J}(t)}{dt}= \int_0^\ell m(x) u^2_{t}dx
- \int_0^\ell \left(r(x)u_{xx}\right)_{xx}u dx \qquad \qquad \qquad \qquad  \nonumber \\ [1pt]
\qquad +k^{-}_a u_{x}(0,t)u_{xt}(0,t)+k^{+}_a u_{x}(\ell,t)u_{xt}(\ell,t),~t\in[0,T].
\end{eqnarray*}
To transform the second right-hand side integral here, we employ the identity
\begin{eqnarray*}
-\int_0^\ell  \left (r(x) u_{xx}\right )_{xx} u dx=
-\int_0^\ell  r(x) u^2_{xx} dx -k^{-}_r u^2_{x}(0,t)-k^{-}_a u_{x}(0,t)u_{xt}(0,t) \\ [1pt]
-k^{+}_r u^2_{x}(\ell,t)-k^{-}_a u_{x}(\ell,t)  u_{xt}(\ell,t),~t\in[0,T].
\end{eqnarray*}
Then we get:
\begin{eqnarray*}
\frac{d \mathcal{J}(t)}{dt}= \int_0^\ell m(x) u^2_{t}dx
-\int_0^\ell  r(x) u^2_{xx} dx -k^{-}_r u^2_{x}(0,t)-k^{+}_r u^2_{x}(\ell,t),
\end{eqnarray*}
for all $t\in[0,T]$. This leads to the required result (\ref{10}).
 \hfill$\Box$

The next lemma shows another relationship between the auxiliary function $\mathcal{J}(t)$ and the energy function $\mathcal{E}(t)$. Namely, it shows that the energy function serves as lower and upper bounds to the auxiliary function introduced in (\ref{9}).

\medskip
\begin{lemma}\label{Lemma-3}
Assume that in addition to the basic conditions (\ref{2}), the coefficient $r(x)$
in (\ref{1}) satisfies the regularity condition: $r \in H^2(0,\ell)$. Then the
following inequalities hold:
\begin{eqnarray}\label{11}
-\beta_0 \,\mathcal{E}(t) \le \mathcal{J}(t) \le \beta_1\, \mathcal{E}(t), ~t\in[0,T],
\end{eqnarray}
where
\begin{eqnarray}\label{12}
\left. \begin{array}{ll}
\displaystyle \beta_0 =\frac{\ell^2}{2}\,\sqrt{\frac{m_1}{r_0}}\,\\ [14pt]
\displaystyle \beta_1=\beta_0 \left \{1+  \frac{1}{\sqrt{m_1 r_0}}\left [\ell^2\mu_1+\frac{2}{\ell} \left (k_a^{-}+k_a^{+} \right )\right ]\right\}\,, \\
\end{array} \right.
\end{eqnarray}
and $m_1,\,\mu_1,\,r_0>0$ are the constants introduced in (\ref{2}).
\end{lemma}
{\bf Proof.} We estimate separately each term on the right hand side of formula (\ref{9}). For the first term we use the $\varepsilon$-inequality to get
\begin{eqnarray}\label{13}
\left \vert \int_0^\ell m(x) u u_{t}dx\right \vert \le \frac{\varepsilon}{2}\, \int_0^\ell m(x) u_t^2dx + \frac{1}{2\varepsilon}\, \int_0^\ell m(x) u^2dx.
\end{eqnarray}
Under the condition $r \in H^2(0,\ell)$ the exists the regular weak solution $u\in L^2(0,T; H^4(0,\ell))$, with $u_t\in L^2(0,T; \mathcal{V}^2(0,\ell))$, $u_{tt}\in L^2(0,T;L^2(0,\ell))$ and $u_{ttt}\in L^2(0,T;H^{-2}(0,\ell))$ of problem (\ref{1}) \cite{Hasanov-Romanov:2021}. For this solution we employ the inequality
\begin{eqnarray}\label{14}
\int_0^\ell u^2 dx \le \frac{\ell^4}{4} \int_0^\ell u_{xx}^2dx,~t\in [0,T],
\end{eqnarray}
which can be easily proved due to the conditions $u(0,t)=u(\ell,t)=0$. This yeilds:
\begin{eqnarray*}
\int_0^\ell m(x) u^2dx \le \frac{\ell^4 \rho_1}{4 r_0} \int_0^\ell r(x)u_{xx}^2dx,
\end{eqnarray*}
Substituting this in (\ref{13}) we get:
\begin{eqnarray*}
\left \vert \int_0^\ell m(x) u u_{t}dx\right \vert \le \frac{\varepsilon}{2}\, \int_0^\ell m(x) u_t^2dx + \frac{\ell^4 m_1}{8\varepsilon r_0} \int_0^\ell r(x)u_{xx}^2dx.
\end{eqnarray*}
Choose here the parameter $\varepsilon>0$ from the condition $\varepsilon/2=\ell^4 m_1/(8 r_0\,\varepsilon)$ as
\begin{eqnarray*}
\varepsilon= \frac{\ell^2}{2}\,\sqrt{\frac{m_1}{r_0}}\,,
\end{eqnarray*}
we obtain the following estimate:
\begin{eqnarray}\label{15}
\left \vert \int_0^\ell m(x) u u_{t}dx\right \vert \le
\frac{\ell^2}{4}\,\sqrt{\frac{m_1}{r_0}}\left [ \int_0^\ell m(x) u_t^2dx + \int_0^\ell r(x) u_{xx}^2dx \right ].
\end{eqnarray}

Now, we estimate the second right hand side integral in formula (\ref{9}), using
inequality (\ref{14}). We have:
\begin{eqnarray}\label{16}
\int_0^\ell \mu(x) u^2dx \le \frac{\ell^4 \mu_1}{4 r_0} \int_0^\ell r(x)u_{xx}^2dx.
\end{eqnarray}

Finally, to estimate the third and fourth terms on the right side of formula (\ref{9}), we use the same argument as above to conclude that
\begin{eqnarray*}
u^2_x(0,t)=\left(-\int_0^{\tilde{x}} u_{xx}(x,t)dx \right)^2\le \tilde{x}\, \int_0^{\tilde{x}} u^2_{xx}(x,t)dx,\\
u^2_x(\ell,t)=\left(\int_{\tilde{x}}^{\ell} u_{xx}(x,t)dx \right)^2\le (\ell-\tilde{x})\, \int_0^{\tilde{x}} u^2_{xx}(x,t)dx.
\end{eqnarray*}
Hence,
\begin{eqnarray}\label{17}
\left. \begin{array}{ll}
\displaystyle \frac{1}{2}\, k^{-}_a  u_x^2 (0,t) \le \frac{\ell}{2}\,\frac{k^{-}_a }{r_0}\int_0^{\ell} r(x)u^2_{xx}(x,t)dx,\\ [9pt]
\displaystyle \frac{1}{2}\, k^{+}_a  u_x^2 (\ell,t)\le \frac{\ell}{2}\,
\frac{k^{+}_a }{r_0}\int_0^{\ell} r(x) u^2_{xx}(x,t)dx.
\end{array} \right.
\end{eqnarray}
In view of (\ref{15}), (\ref{16}) and  (\ref{17}) we obtain the following upper estimate for the auxiliary function $\mathcal{J}(t)$:
\begin{eqnarray*}
\mathcal{J}(t)\le  \frac{\ell^2}{4}\,\sqrt{\frac{m_1}{r_0}} \int_0^\ell m(x) u_t^2dx \qquad \qquad \qquad \qquad \qquad \qquad \qquad \\ [1pt]
\qquad \qquad \qquad  +\left [\frac{\ell^2}{4}\,\sqrt{\frac{m_1}{r_0}}+\frac{\ell^4}{4r_0}\,\mu_1+
\frac{\ell}{2r_0} \left (k^{-}_a+k^{+}_a \right ) \right ]\int_0^\ell r(x) u_{xx}^2dx,
\end{eqnarray*}
for all $t\in (0,T]$. This leads to the upper bound
\begin{eqnarray*}
\mathcal{J}(t) \le \beta_1\, \mathcal{E}(t), ~t\in[0,T],
\end{eqnarray*}
in terms of the energy functional $ \mathcal{E}(t)$ and the constant $\beta_1>0$ introduced in  (\ref{12}).

The lower bound
\begin{eqnarray*}
\mathcal{J}(t) \ge -\beta_0\, \mathcal{E}(t), ~t\in[0,T]
\end{eqnarray*}
follows from the second part
\begin{eqnarray*}
\int_0^\ell m(x) u u_{t}dx \ge -\,\frac{\ell^2}{4}\,\sqrt{\frac{m_1}{r_0}}\left [ \int_0^\ell m(x) u_t^2dx + \int_0^\ell r(x) u_{xx}^2dx \right ]
\end{eqnarray*}
of estimate (\ref{15}).  This leads to the required estimates (\ref{11}). \hfill$\Box$

\medskip
\begin{remark}\label{Remark-3}
The constants $\beta_0,\beta_1>0$ introduced in (\ref{12}) depend only on the geometric and physical parameters of a beam.
\end{remark}

We introduce now the Lyapunov function
\begin{eqnarray}\label{18}
\mathcal{L}(t)=\mathcal{E}(t)+\lambda\mathcal{J}(t),\, t\in[0,T]
\end{eqnarray}
through the energy function $\mathcal{E}(t)$ and the auxiliary function $\mathcal{J}(t)$, where  $\lambda>0$ is the penalty term.

\medskip
\begin{theorem}\label{Theorem-2}
Assume that the inputs in (\ref{1}) satisfy the basic conditions (\ref{2}) and
the regularity condition $r \in H^2(0,\ell)$. Suppose, in addition that the damping constant of proportionality is positive,
\begin{eqnarray}\label{19}
\gamma_0>0.
\end{eqnarray}
Then system (\ref{1}) is exponentially stable, that is,
\begin{eqnarray}\label{20}
\mathcal{E}(t)\le M\,  e^{-\sigma \,t}\, \mathcal{E}(0),~t\in[0,T],
\end{eqnarray}
where
\begin{eqnarray}\label{21}
\left. \begin{array}{ll}
\displaystyle M= \frac{1+  \beta_1 \lambda}{1-  \beta_0 \lambda }~, ~
 \sigma=\frac{2 \lambda}{1+\beta_1 \lambda}~,\\ [14pt]
 0<\lambda <\min (1/ \beta_0,\, \gamma\,m_0/(2m_1)),
 \end{array} \right.
\end{eqnarray}
where $\mu_0,m_1>0$ and $\beta_0, \beta_1>0$ are the constants introduced in (\ref{2}) and (\ref{12}), respectively, and $\mathcal{E}(0)>0$ is the initial energy defined in (\ref{5}).
\end{theorem}
{\bf Proof.} Using estimates (\ref{11}) in (\ref{18}) we get:
\begin{eqnarray*}
\left (1-\beta_0 \lambda \right ) \,\mathcal{E}(t) \le \mathcal{L}(t) \le
\left (1+\beta_1 \lambda \right )\,\mathcal{E}(t), ~t\in[0,T].
\end{eqnarray*}
From the positivity requirement of the first left hand side multiplier, we find that the penalty term should satisfiy the following condition:
\begin{eqnarray}\label{22}
0<\lambda <1/ \beta_0,~\beta_0>0.
\end{eqnarray}
Differentiate now $\mathcal{L}(t)$ with respect to the variable $t\in (0,T)$ and use formulas (\ref{7}) and (\ref{10}). We have:
\begin{eqnarray}\label{23}
\frac{d \mathcal{L}(t)}{dt}+2 \lambda \mathcal{E}(t)= -\int_0^\ell \left [\mu(x)-2\lambda m(x)\right ] u_t^2dx
 \qquad \qquad \quad \nonumber\\ [2pt]
\qquad \qquad \qquad -k^{-}_a u^2_{xt}(\ell,t)-k^{+}_a u^2_{xt}(\ell,t),~t\in[0,T].
\end{eqnarray}
Assume that, in addition to (\ref{22}), the penalty term satisfies also the following condition:
\begin{eqnarray*}
\lambda \le \mu_0/(2m_1)
\end{eqnarray*}
which guarantees positivity of the term in the square bracket under the right hand side intagral in (\ref{23}). In view of the relation $\mu_0=\gamma\,m_0$, this condition implies
\begin{eqnarray}\label{24}
\lambda \le \gamma\,m_0/(2m_1).
\end{eqnarray}
This leads to
\begin{eqnarray*}
\frac{d \mathcal{L}(t)}{dt}+2 \lambda \mathcal{E}(t)\le0,~t\in[0,T],
\end{eqnarray*}
or, with $\mathcal{E}(t) \ge \mathcal{L}(t)/(1+\lambda \gamma_1)$, to the inequality
\begin{eqnarray*}
\frac{d \mathcal{L}(t)}{dt}+\frac{2\lambda}{1+\lambda \gamma_1}\,\mathcal{L}(t)\le 0,~t\in[0,T].
\end{eqnarray*}
Solving this inequality we find:
\begin{eqnarray*}
\mathcal{L}(t) \le e^{-\sigma \,t}\, \mathcal{E}(0),~t\in[0,T]
\end{eqnarray*}
which implies the required estimate (\ref{20}).
\hfill$\Box$

\medskip
\begin{remark}\label{Remark-4}
The constant $\sigma>0$ in (\ref{21}), called the decay rate parameter, depends only on the geometric and physical parameters of the beam and also on the stiffness of the torsional dampers introduced in (\ref{2}), as formulas (\ref{12}) show. Hence, the uniform exponential stability estimate (\ref{21}) can be applied to study exponential stability for Euler-Bernoulli beams with various physical and geometric properties, under boundary controls in rotation and angular velocity. Furthermore, considering formula (\ref{12}), estimate (\ref{21}) also clearly shows the contribution of each damping factor $\mu(x)$, $k_a^{-}$ and $k_a^{-}$ to the energy decay rate.
\end{remark}

\section{Numerical results}

Although there is an exponential function $e^{-\sigma \,t}$ on the right side of the estimate (\ref{20}), with the decay rate parameter $\sigma >0$ introduced in (\ref{21}), in some  cases, this appearance can be misleading. Namely, $\sigma >0$ is dependent on the positive parameters  $\lambda$ and $\beta_1$. The specific values of these parameters play a crucial role in determining the decay behavior of the function $e^{-\sigma \,t}$. Depending on the values of $\lambda$ and $\beta_1$, the decay of this function can exhibit characteristics similar to the decay of a linear function. To see such cases, it is necessary to study the dependence of the decay rate parameter on not only the geometric and physical parameters of the beam, but also on the viscous external damping parameter $\mu(x)$ and the torsional dampers $k_a^{-}, k_a^{-}$ separately.

The examples below are provided to illustrate these situations and their causes. Without loss of generality, here we consider the constant coefficient beam equation
\begin{eqnarray}\label{25}
m u_{tt}+\mu u_{t}+r u_{xxxx}=0,\, (x,t) \in \Omega_T,
\end{eqnarray}
where 
\begin{eqnarray*}
m= \rho\, S, ~\mu=\gamma m, ~r= E I,
\end{eqnarray*}
in accordance with the above notation. For this constant coefficients equation, formulas (\ref{12}) and (\ref{21}) for the parameters $\beta_0,\, \beta_1,\,M_1,\,\sigma>0$ and conditions are as follow:
\begin{eqnarray}\label{26}
\left. \begin{array}{ll}
\displaystyle \beta_0 =\frac{\ell^2}{2}\,\sqrt{\frac{m}{r}}\,\\ [14pt]
\displaystyle \beta_1=\beta_0 \left [1+  \ell^2 \sqrt{\frac{m}{r}}\,\gamma\right ]
+\frac{\ell}{2r} \left (k_a^{-}+k_a^{+} \right )\,, \\ [14pt]
\displaystyle M= \frac{1+  \beta_1 \lambda}{1-  \beta_0 \lambda }~, ~
 \sigma=\frac{2 \lambda}{1+\beta_1 \lambda}~.
\end{array} \right.
\end{eqnarray}

Here, the beam with the rectangular cross section $S=b\, h$, where $b>0$ and $h>0$ are the width height, with the following numerical values of the geometric and physical parameters are examined \cite{Repetto:2012}:
\begin{eqnarray}\label{27}
\left. \begin{array}{ll}
\ell=0.502\,\mbox{m},~b=1.7\times 10^{-3}\,\mbox{m},~h=0.89\times 10^{-3} \,\mbox{m},\\ [4pt]
\rho=1.42\times 10^{3}\,\mbox{Kg\,m}^{-3},~E=3.1\times 10^{9}\,\mbox{N/m}^{2},~\gamma \in [0.01,\,10]\, \mbox{s}^{-1}.
 \end{array} \right.
\end{eqnarray}

With the numerical values in (\ref{27}) we have:
\begin{eqnarray*}
\left. \begin{array}{ll}
S=1.51\times 10^{-6}\,\mbox{m}^2,~I:=bh^3/12=0.1\times 10^{-12}\,\mbox{m}^3,\\ [4pt]
m=2.14\times 10^{-3}\,\mbox{Kg\,m}^{-1},~r=0.31\times 10^{-3} \,\mbox{N\,m}^2,~\mu=0.22\,\mbox{Kg\,m}^{-1}.
 \end{array} \right.
\end{eqnarray*}

We consider three-level, weak, medium, and high damping cases corresponding to the values $\gamma=0.1$, $\gamma=1.0$ and $\gamma=5.0$ of the damping constant of proportionality,
using the following values $\langle k^{-}_a,k^{+}_a \rangle= \langle 0,\, 0 \rangle$ and  $\langle k^{-}_a,k^{+}_a \rangle= \langle 0.01,\, 0.01 \rangle$ of the stiffness of the torsional dampers.

The calculated by formulas given in (\ref{26}) values of the decay rate parameter $\sigma>0$ are listed in Table 1. The values of the penalty term $\lambda>0$ are set according to the requirement $0<\lambda <\min (1/ \beta_0,\, \gamma/2)$.

From the last column of Table 1 it can be seen that, in absence of the torsional dampers ($k^{-}_a=k^{+}_a=0$), the increase in the value of the damping constant from $\gamma=0.1$ to $\gamma=5.0$, leads to the increase of the decay parameter $\sigma>0$. Thus, for the
weak damping case $\gamma=0.01$ the value of the decay parameter is $\sigma=0.08$, and
the energy decay is only exponential in appearance, in fact, it is linear (Figure 1 on the left).

\begin{center}\label{Table-1}
{\textbf{Table 1.} The decay rate parameters corresponding to the geometric and physical parameters given in (\ref{27}).}\\
\vskip .2in
\begin{tabular}{|c|}
\hline Damping constant \\
\hline  \\
$\gamma=0.1$  \\
\hline  \\
$\gamma=1.0$  \\
 \hline  \\
$\gamma=5.0$
 \\ \hline
\end{tabular}
\begin{tabular}{|c|c|c|c|c|}
\hline  $\langle k^{-}_a,k^{+}_a \rangle$ & $\langle \beta_0,\beta_1\rangle$   & $\lambda$  & $M$ & $\sigma$      \\ \hline
$\langle 0,\,0 \rangle$ & $\langle 0.33,\,0.35 \rangle$  & $0.04$ & $1.03$ & $0.08$    \\ \hline
 $\langle 0.01,\,0.01\rangle$ & $\langle 0.33,\,16.55 \rangle$ & $0.04$ & $1.68$ & $0.05$  \\ \hline
\hline  $\langle 0,\,0 \rangle$ & $\langle 0.33,\,0.55 \rangle$ & $0.4$ & $1.41$ & $0.66$    \\ \hline
 $\langle 0.01,\,0.01 \rangle$ & $\langle 0.33,\,16.75 \rangle$ & $0.4$ & $8.87$ & $0.10$  \\ \hline
\hline $\langle 0,\,0 \rangle$ & $\langle 0.33,\,1.42 \rangle$ & $2.4$ & $21.19$ & $1.09$    \\ \hline
$\langle 0.01,\,0.01 \rangle$ & $\langle 0.33,\,17.62  \rangle$  & $2.4$ & $208.12$ & $0.11$  \\ \hline
\end{tabular}
\end{center}

\begin{figure}[!htb]
\includegraphics[width=7.0cm,height=7.0cm]{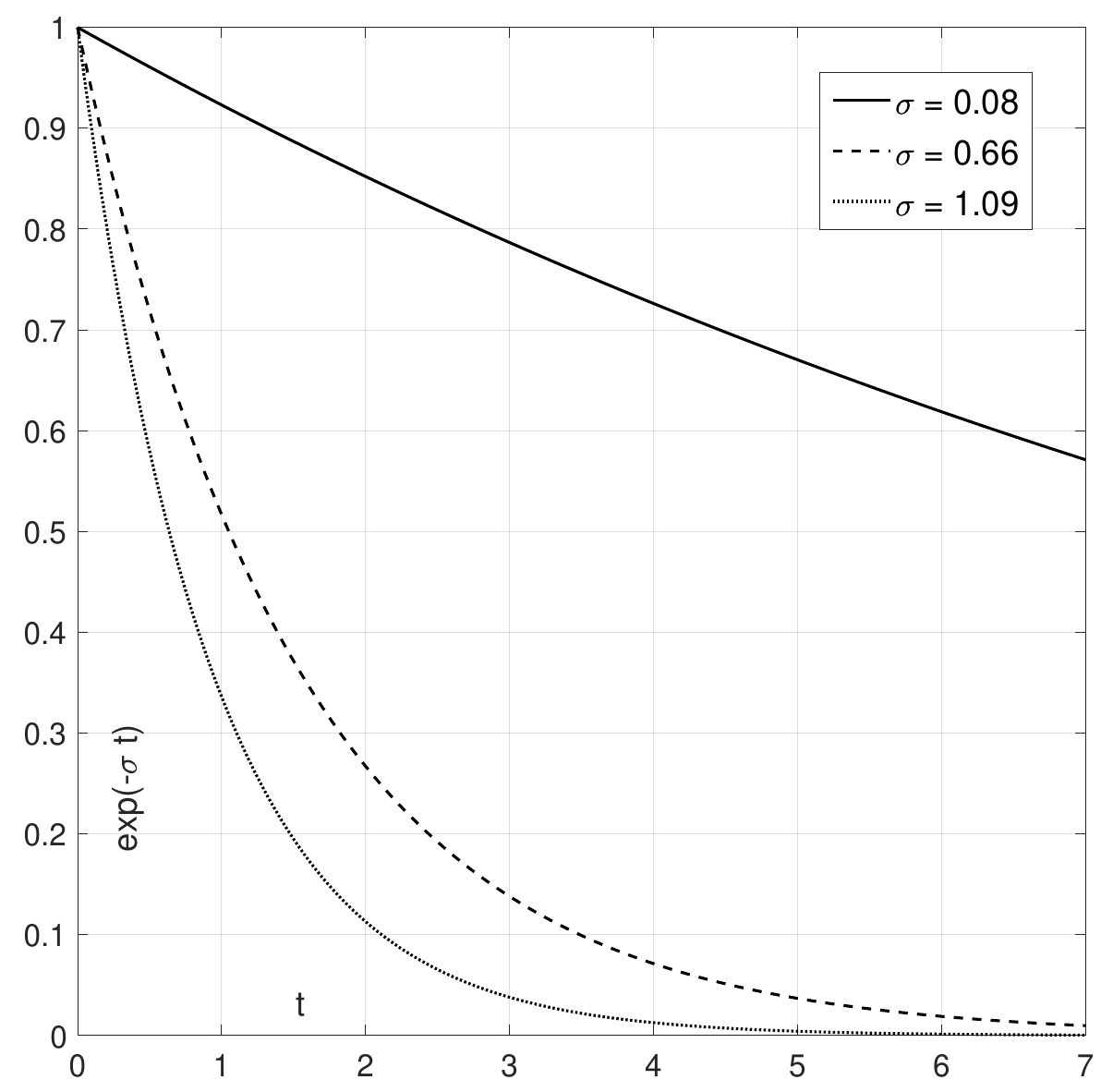}~~
\includegraphics[width=7.0cm,height=7.0cm]{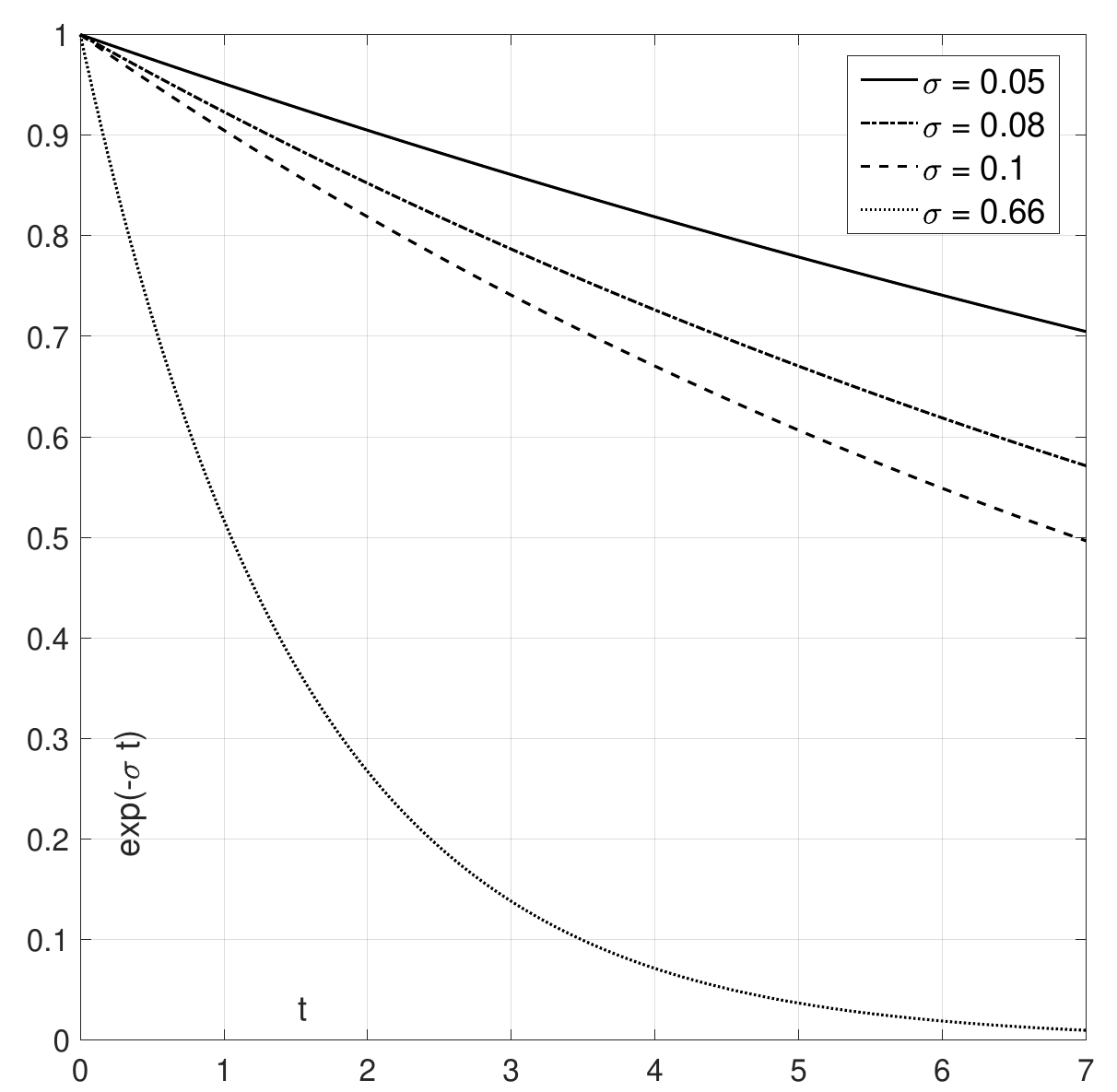}
\noindent \caption{Behaviour of the function $\exp(-\sigma t)$: with 
$k^{-}_a=k^{+}_a =0$ (left) and with with
$k^{-}_a=k^{+}_a =0.01$ (right).} \label{Fig1ab}
\end{figure}

Comparing the values of the decay rate parameter, in the last column of Table 1, corresponding to zero and non-zero values of the stiffness of the torsional dampers, we can observe the role of these boundary controls (Figure 1 on the right).

\section{Conclusions}
 This study proposes an approach for the exponential stability analysis of Euler-Bernoulli beams under boundary controls in rotation and angular velocity. By employing the regular weak solution, energy identity, and Lyapunov function, we are able to derive a uniform exponential decay estimate for the system's energy.
 
Our approach is grounded in natural assumptions concerning physical parameters and other inputs, ensuring the existence of a regular weak solution. The decay rate constant in the derived estimate relies solely on the physical and geometric parameters of the beam, which include the viscous external damping coefficient, as well as the boundary springs and dampers. This feature enables straightforward utilization of decay rate estimation in practical engineering applications.

Furthermore, we have provided preliminary numerical examples that shed light on the role of damping parameters. However, a more detailed analysis, focusing on the individual contributions of each damping parameter to the overall damping behavior, will be pursued in future research.

\section*{Acknowledgments}
The research has been supported by the Scientific and Technological Research Council of Turkey (TUBITAK) through the Incentive Program for International Scientific Publications (UBYT). The research of the author has also been supported by FAPESP, through the Visiting Researcher Program, proc. 2021/08936-1, in Escola Polit\'{e}cnica, University of S\~{a}o Paulo, Brazil, during the period November 02 - December 18, 2022.



\begin{thebibliography}{16}
\bibitem{Banks:Inman:1991}
H.T. Banks, D.J. Inman, On Damping Mechanisms in Beams, Journal of Applied Mechanics, 58(3) (1991) 716--723.
\bibitem{Cai:2022}
J. Cai, P. Le Grognec, Lateral buckling of submarine pipelines under high temperature and high pressure - A literature review, Ocean Engineering 244(15) (2022) 110254.
\bibitem{Chen-Krantz:1988}
G. Chen, S.G. Krantz, D.W. Ma, C.E. Wayne, H.H. West, H. The Euler-Bernoulli beam equation with boundary energy dissipation. Report, 1 Sep. 1985 - 31 Aug. 1987, Pennsylvania State Univ., University Park., 1, 1988. https://dx.doi.org/10.21236/ada189517.
\bibitem{Chen-Russell:1982}
G. Chen, D.L. Russell, A mathematical model for linear elastic systems with structure
damping, Quart. Appl. Math. 39(1982) 433--454.
\bibitem{Chen-Xu:2014}
Y.L. Chen, G. Q. Xu, Exponential stability of uniform Euler-Bernoulli beams with non-collocated boundary controllers, J. Math. Anal. Appl. 409(2014) 851--867.
\bibitem{Clough-Penzien:1975}
R.C. Clough, J. Penzien, Dynamics of Structures, McGraw Hill Inc., New York, 1975.
\bibitem{Crandall:1970}
S.H. Crandall, The Role of Damping in Vibration Theory, J. Sound Vibr. 11(1970) 3--18.
1970.
\bibitem{Evans:2002}
L.C. Evans, Partial Differential Equations, 2nd edn., American Mathematical Society, Rhode Island, 2010.
\bibitem{Grognec:2020}
P. Le Grognec, A. N\'{e}me, J. Cai, Investigation of the torsional effects on the lateral buckling of a pipe-like beam resting on the ground under axial compression, International Journal of Structural Stability and Dynamics 20 (9) (2020) 2050110.
\bibitem{Guo:2001}
B.-Z. Guo and R. Yu, On Riesz basis property of discrete operators with application to an
Euler-Bernoulli beam equation with boundary linear feedback control, IMA J. Math. Control
Inform. 18 (2001) 241--251.
\bibitem{Guo:2002}
B.-Z. Guo, Riesz basis property and exponential stability of controlled Euler–Bernoulli beam
equations with variable coefficients, SIAM J. Control Optim. 40 (2002) 1905--1923.
\bibitem{F-Guo:2004}
F Guo, F Huang, Boundary Feedback Stabilization of the Undamped Euler--Bernoulli Beam with Both Ends Free, SIAM J. Control Optim. 43(1) (2004) 341--356.
\bibitem{Hasanov-Romanov:2021}
A. Hasanov Hasanoglu, A.G. Romanov, Introduction to Inverse Problems for Differential Equations, 2nd ed, Springer, New York, 2021.
\bibitem{Hong:2015}
Z. Hong, R. Liu, W. Liu, S. Yan, A lateral global buckling failure envelope for a high temperature and high pressure (ht/hp) submarine pipeline, Applied Ocean Research 51 (2015)  117–128.
\bibitem{Huang:1986}
F.L. Huang, Some problems for linear elastic systems with damping, Acta Math. Sci. 6
(1986) 101--107.
\bibitem{Huang:1988}
F. Huang, On the mathematical model for linear elastic systems with analytic damping, SIAM
J. Control Optim. 26 (1988) 714--724.
\bibitem{Inman:2015}
D. J. Inman, Engineering Vibration, 4th Edn., Pearson Education Limited, 2014.
\bibitem{Liu:2018}
R. Liu, X. Wang, Lateral global buckling high-order mode analysis of a submarine pipeline with imperfection, Applied Ocean Research 73 (2018) 107--126.
\bibitem{Park:2019}
Y.S. Park, S. Kim, N. Kim, J.J. Lee, Evaluation of bridge support condition using bridge responses. Structural Health Monitoring, 18(3) (2019) 767-777.
\bibitem{Repetto:2012}
C.E. Repetto, A. Roatta and R.J. Welti, Forced vibrations of a cantilever beam, Eur. J. Phys. 33 (2012) 1187--1195.
\bibitem{Russell:1978}
D. L. Russell, Controllability and stabilizatiblity theory for linear partial differential equations: Recent progress and open questions, SIAM Rev., 20 (1978) 639--739.

\end{thebibliography}
\end{document}